\documentclass[reqno,a4paper, 12pt]{amsart}

\usepackage[pdfpagelabels]{hyperref}
\usepackage{graphicx}

\usepackage[ansinew]{inputenc}
\usepackage{amsfonts,epsfig}
\usepackage{latexsym}
\usepackage{amsmath}
\usepackage{amssymb}

\newtheorem{theorem}{Theorem}

\newtheorem{corollary}[theorem]{Corollary}

\newtheorem{lettertheorem}{Theorem}
\newtheorem{letterlemma}[lettertheorem]{Lemma}

\theoremstyle{definition}

\theoremstyle{remark}

\newtheorem{example}{Example}

\numberwithin{equation}{section}

\newcommand{\set}[1]{\left\{#1\right\}}
\newcommand{\abs}[1]{\lvert#1\rvert}
\newcommand{\nm}[1]{\lVert#1\rVert}

\newcommand{\D}{\mathbb{D}}

\newcommand{\N}{\mathbb{N}}

\newcommand{\C}{\mathbb{C}}

\newcommand{\e}{\varepsilon}

\renewcommand{\phi}{\varphi}

\renewcommand{\rho}{\varrho}
\DeclareMathOperator{\Real}{Re}

\def\VMOA{\mathord{\rm VMOA}}

\def\d{\delta}           \def\e{\varepsilon}

         \def\r{\rho}         \def\z{\zeta}
                  \def\vp{\varphi}

\begin{document}

\title[On oscillation of solutions]{On oscillation of solutions of linear differential equations}
\thanks{The first author is supported by the Academy of Finland \#258125,  and the second
author is supported in part by the Academy of Finland \#268009, the 
Faculty of Science and Forestry of University of Eastern Finland \#930349, and the grants 
MTM2011-26538 and MTM2014-52865-P}

\author{Janne Gr\"ohn}
\address{Department of Physics and Mathematics, University of Eastern Finland, P.O. Box 111, FI-80101 Joensuu, Finland}
\email{janne.grohn@uef.fi}

\author{Jouni R\"atty\"a}
\address{Department of Physics and Mathematics, University of Eastern Finland, P.O. Box 111, FI-80101 Joensuu, Finland}
\email{jouni.rattya@uef.fi}



\date{\today}


\begin{abstract}
An interrelationship is found
between the accumulation points of zeros of non-trivial solutions of $f''+Af=0$ and the boundary
behavior of the analytic coefficient~$A$ in the unit disc $\D$ of the complex plane $\C$.

It is also shown that the geometric distribution of zeros of any non-trivial solution of $f''+Af=0$ 
is severely restricted if 
\begin{equation} \label{eq:cs_a}\tag{$\star$}
|A(z)| (1-|z|^2)^2 \leq 1 + C (1-|z|), \quad z\in\D,
\end{equation}
for any constant $0<C<\infty$. These considerations are related to the open problem
whether \eqref{eq:cs_a} implies finite oscillation for all non-trivial solutions.
\end{abstract}

\maketitle


\vspace*{0.3cm}


\section{Introduction}

The following result plays a decisive role in the oscillation theory
of solutions of linear differential equation
\begin{equation} \label{eq:de2}
    f'' + A f = 0
    \end{equation}
in the unit disc $\D$ of the complex plane $\C$.
If $A$ is an~analytic function in $\D$ for which
    \begin{equation}\label{Eq:Nehari}
      |A(z)| (1-|z|^2)^2 \leq 1, \quad z\in\D,
    \end{equation}
then each non-trivial solution $f$ of \eqref{eq:de2} vanishes at most once in $\D$. This
statement corresponds to the well-known result of
Z.~Nehari~\cite[Theorem~1]{N:1949}, which provides a sufficient condition for
injectivity of any locally univalent meromorphic function $w$ in $\D$ in terms of the size of its Schwarzian derivative
\begin{equation*}
  S_w=\frac{w'''}{w'}-\frac32\left(\frac{w''}{w'}\right)^2=\left(\frac{w''}{w'}\right)'-\frac12\left(\frac{w''}{w'}\right)^2.
\end{equation*}

The corresponding necessary condition was invented by
W.~Kraus~\cite{Kraus1932}, and rediscovered by
Nehari~\cite[Theorem~1]{N:1949} some years later. In the setting of differential
equations it states that, if $A$ is analytic in $\D$, and each solution $f$ of \eqref{eq:de2}
vanishes at most once in $\D$, then $|A(z)| (1-|z|^2)^2 \leq 3$ for all $z\in\D$. An important 
discovery of B.~Schwarz \cite[Theorems~3--4]{S:1955} shows that the condition 
\begin{equation*} 
\sup_{z\in\D}|A(z)|(1-|z|^2)^2<\infty,
\end{equation*}
which allows non-trivial solutions of \eqref{eq:de2} to have infinitely many zeros in $\D$,
is both necessary and sufficient for zeros of all non-trivial solutions to be 
separated with respect to the hyperbolic metric. 

Our first objective is to consider the interrelationship
between the accumulation points of zeros of non-trivial solutions $f$ of \eqref{eq:de2} and the boundary
behavior of the coefficient $A$.
The second objective is a question of more specific nature.
We consider differential equations \eqref{eq:de2} in which
the growth of the coefficient barely exceeds the bound \eqref{Eq:Nehari} that ensures finite oscillation.


\section{Results}


\subsection{Accumulation points of zeros of solutions} \label{subsec:geo}

The point of departure is a~result, which associates
the zero-sequences of non-trivial solutions of \eqref{eq:de2} to the boundary behavior of the
coefficient. This theorem sets the stage for more profound oscillation theory.


\begin{theorem}\label{Thm-Nehari-Local}
Let $A$ be an analytic function in $\D$, and let $\z\in\partial\D$.

If there exists a sequence $\{w_n\} \subset \D$ converging to $\zeta$, such that
\begin{equation}\label{Eq-Local-Nehari-Improved}
  |A(w_n)|(1-|w_n|^2)^2 \to c
\end{equation}
for some $c\in (3,\infty]$, then for each $\delta> 0$ there exists a
non-trivial solution of \eqref{eq:de2} having two distinct zeros in $D(\zeta,\delta) \cap \D$.

Conversely, if for each $\delta>0$ there exists a
non-trivial solution of \eqref{eq:de2} having two distinct zeros in $D(\zeta,\delta) \cap \D$,
then there exists a sequence $\{w_n\} \subset \D$ converging to $\zeta$ such that
\eqref{Eq-Local-Nehari-Improved} holds for some $c\in [1,\infty]$.
\end{theorem}

We point out that \eqref{Eq-Local-Nehari-Improved} with $c\in (3,\infty]$ does not necessarily imply
infinite oscillation for any non-trivial solution of \eqref{eq:de2}, see Example~\ref{ex:zeropairs} below.

The proof of Theorem~\ref{Thm-Nehari-Local} is based on theorems by Nehari and
Kraus, and on a~principle of localization. One of the key factors is an application of a suitable
family of conformal maps under which the image of $\D$ has 
a smooth boundary,  that intersects $\partial\D$ precisely on an
arc centered at $\z\in\partial\D$. The second assertion of Theorem~\ref{Thm-Nehari-Local}
is implicit in the proof of \cite[Theorem~1]{N:1949}, and follows directly from the following property:
if $z_1,z_2\in \D$ are two distinct zeros of a non-trivial solution $f$
of \eqref{eq:de2}, then 
there exists a point $w\in\D$, which belongs to the hyperbolic geodesic going through $z_1$ and $z_2$, such
that $\abs{A(w)} (1-\abs{w}^2)^2 > 1$.


\subsection{Chuaqui-Stowe question} \label{subsec:cs}

Schwarz~\cite[Theorem~1]{S:1955} supplemented the oscillation theory by
proving that, if there exists a constant $0<R<1$ such that 
\begin{equation*}
|A(z)| (1-|z|^2)^2 \leq 1, \quad R<|z|<1,
\end{equation*} 
then each non-trivial solution of \eqref{eq:de2} has at most finitely many zeros.
Schwarz also gave an example \cite[p.~162]{S:1955} showing that
the constant one in the right-hand side of \eqref{Eq:Nehari} is best possible.
That is,
for each $\gamma>0$, the functions
\begin{equation*}
A(z) = \frac{1 + 4\gamma^2}{(1-z^2)^2}\quad \text{and} \quad
f(z) = \sqrt{1-z^2} \, \sin \!\left( \gamma \log\frac{1+z}{1-z} \right)
\end{equation*}
satisfy \eqref{eq:de2}, while $f$ has infinitely many (real) zeros in $\D$.
Example~\ref{ex:KrausSharpness} below shows the sharpness of Kraus' result.

M.~Chuaqui and D.~Stowe~\cite[Theorem~5]{CS:2008}
constructed an example showing that for each continuous function
$\varepsilon\colon [0,1) \to [0,\infty)$ satisfying
$\varepsilon(r) \to \infty$ as $r\to 1^-$, there exists an analytic function $A$ such that
    \begin{equation}\label{eq:bound}
    \abs{A(z)}  (1-\abs{z}^2)^2  \leq 1+ \varepsilon(\abs{z})  (1-\abs{z}),\quad
    z\in\D,
    \end{equation}
while \eqref{eq:de2} admits a non-trivial solution having infinitely many zeros.
In other words, if $\varepsilon(r) (1-r)$ in \eqref{eq:bound} does not decay to zero
as fast as linear rate as $r\to 1^-$, then non-trivial solutions of \eqref{eq:de2} may
have infinitely many zeros. This is in contrast to the case of real
differential equations \eqref{eq:de2} on the open interval $(-1,1)$,
since then $\varepsilon(r)=(1-r)^{-1} (-\log (1-r))^{-2}$
distinguishes finite and infinite oscillation, see \cite{CDOS:2009,CS:2008} for more details.
Chuaqui and Stowe~\cite[p.~564]{CS:2008} left open a question whether
    \begin{equation} \label{eq:cs}
    \abs{A(z)}  (1-\abs{z}^2)^2 \leq 1 + C (1-\abs{z}), \quad z\in\D,
    \end{equation}
with some or any $0<C<\infty$, implies finite oscillation for all non-trivial solutions
of \eqref{eq:de2}. The following results do not give a complete answer to 
this question, however, they indicate that 
both the growth and the zero distribution of non-trivial
solutions of \eqref{eq:de2} are severely restricted if \eqref{eq:cs} holds for some $0<C<\infty$.


\subsubsection{Growth of solutions}

An analytic function $f$ in $\D$ belongs to
the growth space $H^\infty_\alpha$ for $0\leq \alpha<\infty$, if
\begin{equation*}
\nm{f}_{H^\infty_\alpha} = \sup_{z\in\D} \, \abs{f(z)}(1-\abs{z}^2)^\alpha  < \infty.
\end{equation*}
It is known that the growth of $A$ restricts the growth
of solutions of \eqref{eq:de2}. If $A\in H^\infty_2$, then
there exists a~constant $p=p(\nm{A}_{H^\infty_2})$ with $0\leq p<\infty$ such that
all solutions $f$ of \eqref{eq:de2} satisfy $f\in H^\infty_p$.
This result can be deduced by using classical comparison theorems \cite[Example~1]{P:1982}, 
Gronwall's lemma \cite[Theorem~4.2]{H:2000} or successive approximations 
\cite[Theorem~I]{G:2011}, for example.
We conclude this result by means of straightforward integration. See Example~\ref{ex:growth}
for sharpness discussion.


\begin{theorem} \label{thm:growth}
Let $A$ be an analytic in $\D$ such that
$|A(z)| (1-|z|^2)^2 \leq K + \varepsilon(|z|)$ for all $z\in\D$,
where $0\leq K < \infty$ is a constant, and $\varepsilon(|z|)\to 0$ as $|z|\to 1^{-}$. Then all
solutions of \eqref{eq:de2} belong to $H^\infty_p$ for any $(\sqrt{1+K}-1)/2<p<\infty$.
\end{theorem}


\subsubsection{Separation of zeros of solutions}

The following result establishes a connection between the
separation of zeros of non-trivial solutions of \eqref{eq:de2} and the growth
of the coefficient function~$A$; compare to \cite[Theorems~3 and 4]{S:1955}.

If $z_1,z_2$ are two distinct points in $\D$, then the pseudo-hyperbolic distance
$\varrho_p(z_1,z_2)$ and the hyperbolic distance $\varrho_h(z_1,z_2)$ between $z_1$ and $z_2$
are given by
\begin{equation*}
  \varrho_p(z_1,z_2) = \big| \varphi_{z_1}(z_2) \big|, 
  \quad \varrho_h(z_1,z_2) = \frac{1}{2} \log \frac{1+\varrho_p(z_1,z_2)}{1-\varrho_p(z_1,z_2)},
\end{equation*}
where $\varphi_a(z)= (a-z)/(1-\overline{a} z)$, $a\in\D$. Moreover, let $\xi_h(z_1,z_2)$ denote
the hyperbolic midpoint between $z_1$ and $z_2$. Correspondingly,
\begin{equation*}
\Delta_p(a,r) = \big\{ z\in\D : \varrho_p(z,a)<r \big\},
\quad \Delta_h(a,r) = \big\{ z\in\D : \varrho_h(z,a)<r \big\},
\end{equation*}
are the pseudo-hyperbolic and hyperbolic discs of radius $r>0$ centered at $a\in\D$, respectively.


\begin{theorem}\label{thm:separation}
Let $A$ be an analytic function in $\D$. 

If the coefficient $A$ satisfies \eqref{eq:cs} for some $0<C<\infty$, 
then the hyperbolic distance between any distinct
zeros $z_1,z_2\in\D$ of any non-trivial solution of
\eqref{eq:de2}, for which $1-|\xi_h(z_1,z_2)| < 1/C$, satisfies
    \begin{equation} \label{eq:separation}
    \rho_h(z_1,z_2) \geq \log \frac{2 - C^{1/2}(1-|\xi_h(z_1,z_2)|)^{1/2}}{C^{1/2}(1-|\xi_h(z_1,z_2)|)^{1/2}}.
    \end{equation}

Conversely, if there exists a constant $0<C<\infty$ such that any two distinct zeros 
$z_1,z_2\in\D$ of any non-trivial solution of \eqref{eq:de2}, for which $1-|\xi_h(z_1,z_2)| < 1/C$,
satisfies \eqref{eq:separation}, then
\begin{equation} \label{eq:hypconc}
\abs{A(z)}  (1-\abs{z}^2)^2 \leq 3 \left( 1 + \Psi_C(|z|) (1-|z|)^{1/3} \right), \quad 1-|z|<(8C)^{-1},
\end{equation}
where $\Psi_C$ is positive, and satisfies $\Psi_C(|z|)\longrightarrow \big(2\, (8C)^{1/3}\big)^{+}$ as $|z|\to 1^{-}$.
\end{theorem}

Concerning Theorem~\ref{thm:separation} note that, if $1-|\xi_h(z_1,z_2)| < 1/C$, then \eqref{eq:separation} implies
\begin{equation*} 
  \varrho_h(z_1,z_2) \geq \frac{1}{2} \log \frac{1}{C} + \frac{1}{2} \log \frac{1}{1-|\xi_h(z_1,z_2)|},
\end{equation*}
and hence $\varrho_h(z_1,z_2)$ is large whenever $\xi_h(z_1,z_2)$ is close to the boundary $\partial\D$.

If the coefficient $A$ satisfies \eqref{eq:cs} for some $0<C<\infty$, 
and $f_1$ and $f_2$ are linearly independent solutions of \eqref{eq:de2},
then
the quotient $w=f_1/f_2$ is a normal function (in the sense
of Lehto and Virtanen) by \cite[Corollary, p.~328]{S:1983}. 
As a direct consequence we deduce the following corollary, which states 
that the zero-sequences of $f_1$ and $f_2$ are hyperbolically separated from each other,
see also Example~\ref{ex:csfinite} below.


\begin{corollary} \label{cor:st}
Let $A$ be an analytic function in $\D$, which satisfies \eqref{eq:cs} for some $0<C<\infty$,
and let $\{ z_n \}$ and $\{ \zeta_m \}$ be the zero-sequences of two linearly independent solutions
$f_1$ and $f_2$ of \eqref{eq:de2}. Then, there is a constant $\delta=\delta(f_1,f_2)$
such that $\varrho_h(z_n,\zeta_m)>\delta>0$ for all $n$ and $m$.
\end{corollary}

The following result shows that, if \eqref{eq:cs} does not imply finite oscillation
for non-trivial solutions of \eqref{eq:de2},
then infinite zero-sequences tend to $\partial\D$ tangentially.
Any disc $D(\zeta,1-\abs{\zeta})$ for $\zeta\in\D$,
which is internally tangent to $\D$, is called a horodisc.
Note that Theorem~\ref{th:horodiscs} remains valid in the limit case $C=0$
by the classical theorems of Nehari and Kraus.


\begin{theorem} \label{th:horodiscs}
Let $A$ be an analytic function in $\D$.

If $A$ satisfies \eqref{eq:cs} for some $0<C<\infty$, then any non-trivial solution of \eqref{eq:de2}
has at most one zero in any Euclidean disc $D(\zeta,(1+C)^{-1})$ for $\abs{\zeta} \leq C/(1+C)$.

Conversely, if there exists $0<C<\infty$ such that any non-trivial
solution of \eqref{eq:de2} has at most one zero in any Euclidean disc  $D( \zeta ,(1+C)^{-1})$
for $\abs{\zeta} \leq C/(1+C)$, then
    \begin{equation*}
    |A(z)| (1-|z|^2)^2 \leq 
      3 \Big( 1+ \Psi_C(\abs{z}) \,(1-|z|) \Big), \quad   \frac{C}{1+C} < |z| < 1,
    \end{equation*}
where $\Psi_C$ is positive, and satisfies $\Psi_C(|z|)\longrightarrow ( 2C )^{+}$ as $|z|\to 1^-$.
\end{theorem}

\subsubsection{Geometric distribution of zeros of solutions}

The set 
\begin{equation*} 
Q = Q(I) = \big\{ re^{i\theta}  : e^{i\theta}\in I, \, 1-|I|\leq r < 1\big\}
\end{equation*}
is called a Carleson square based on the arc $I\subset \partial\D$, where $|I| = \ell(Q)$
denotes the normalized arc length of $I$ (i.e.,~$|I|$ is the Euclidean arc length of $I$ divided by $2\pi$).


\begin{theorem} \label{thm:carleson}
If $A$ is an analytic function in $\D$ such that \eqref{eq:cs} holds for some $0<C<\infty$,
then the zero-sequence $\{ z_n \}$ of any non-trivial solution of \eqref{eq:de2} satisfies
\begin{equation} \label{eq:carleson}
\sum_{z_n\in Q} (1-|z_n|)^{1/2} \leq K \, \ell(Q)^{1/2},
\end{equation}
for any Carleson square $Q$. Here $K=K(C)$ with $0<K<\infty$ is a constant independent of $f$.
\end{theorem}

If $A$ is analytic in $\D$ and satisfies \eqref{eq:cs} for some $0<C<\infty$,
then the zero-sequence of any non-trivial solution of \eqref{eq:de2} is interpolating by Theorem~\ref{thm:carleson},
because the zero-sequences are separated by the classical result of Schwarz.


\section{Examples} \label{sec:examples}

We turn to consider some non-trivial examples, the first of which
shows the sharpness of Kraus' result.


\begin{example} \label{ex:KrausSharpness}
Let
\begin{equation*}
A(z) = - \frac{3}{4(1-z)^2}, \quad z\in\D.
\end{equation*}
A solution base $\set{f_1,f_2}$ of \eqref{eq:de2} is given by the non-vanishing functions
\begin{equation*}
f_1(z) = (1-z)^{-1/2},
\quad
f_2(z) = (1-z)^{3/2},
\quad z\in\D.
\end{equation*}
Let $f$ be any non-trivial solution of \eqref{eq:de2}. If $f$ is linearly dependent to $f_1$ or $f_2$, then $f$ is non-vanishing.
Otherwise, there exist $\alpha,\beta\in\C\setminus \set{0}$ such that $f(z)=\alpha f_1(z) + \beta f_2(z)$, and
$f(z)=0$ if and only if $(1-z)^2 = -\alpha/\beta$. This equation has two solutions $z_1,z_2\in\C$, and only one
zero of $f$, say $z_1$, satisfies $\Real{z_1}<1$. This follows from the fact that $1-z_1 = z_2 - 1$. Consequently,
each solution of \eqref{eq:de2} has at most one zero in $\D$, while the coefficient function $A$ satisfies
$\abs{A(z)} (1-|z|^2)^2 \to 3$ as $z\to 1^-$ along the positive real axis.
\end{example}

\smallskip

To conclude that \eqref{Eq-Local-Nehari-Improved} with $c\in
(3,\infty]$ does not imply infinite oscillation for any solution
of \eqref{eq:de2}, we recall Hille's example
\cite[Eq.~(2.12)]{S:1955}. The same example is also
used in \cite[Example~20]{CGHR:2012}.


\begin{example} \label{ex:zeropairs}
Let $A(z) = a/(1-z^2)^2$, where $-\infty<a<0$ is a real parameter. If 
\begin{equation*}
f_1(z) = \sqrt{1-z^2} \left( \frac{1-z}{1+z} \right)^{\frac{1}{2}\sqrt{1-a}},
\quad
f_2(z) = \sqrt{1-z^2} \left( \frac{1-z}{1+z} \right)^{-\frac{1}{2}\sqrt{1-a}},
\end{equation*}
then $\set{f_1,f_2}$ is a solution base of \eqref{eq:de2} of non-vanishing functions.
Let $f$ be any non-trivial solution of \eqref{eq:de2}. If $f$ is linearly dependent to $f_1$ or $f_2$, then $f$ is non-vanishing.
Otherwise, there exist $\alpha,\beta\in\C\setminus \set{0}$ such that $f=\alpha f_1 + \beta f_2$. In this case
$f(z)=0$ if and only if
\begin{equation} \label{eq:zerosoff}
\frac{-\beta}{\alpha} = \left( \frac{1-z}{1+z} \right)^{\sqrt{1-a}}.
\end{equation}

Since $z\mapsto (1-z)/(1+z)$ maps $\D$ onto the right half-plane, 
$\sqrt{1-a} \leq 4$ ensures that each solution of \eqref{eq:de2} has at most two zeros
in $\D$, see also \cite[p.~174]{S:1955}. Further, if  $2 < \sqrt{1-a}$, then there
exists a solution having exactly two zeros in $\D$. In particular, if $-\beta/\alpha$ is real
and strictly negative, then $\alpha f_2 + \beta f_2$ has two zeros in $\D$ by \eqref{eq:zerosoff},
and these zeros are complex conjugate numbers in $\D$.
Note that $2 < \sqrt{1-a} \leq 4$ if and only if $-15\leq a < -3$, and then
$\abs{A(x)} (1-\abs{x}^2)^2 = \abs{a} > 3$ for all $x\in(0,1)$.

We fix $a=-8$, and discuss the zeros of the solution $f = f_1 + k f_2$ for $k>0$. By~\eqref{eq:zerosoff}, the zeros of $f$
in $\D$ are solutions of $-k=(1-z)^3/(1+z)^3$. We conclude that $f$ has exactly two zeros in $\D$
given by
\begin{equation*}
z_1 = \frac{1 - \sqrt[3]{k} \, \exp\left(i \, \pi/3 \right)}{1 + \sqrt[3]{k} \, \exp\left( i \, \pi/3 \right)},
\quad
z_2 = \frac{1 - \sqrt[3]{k} \, \exp\left(-i \, \pi/3\right) }{1 + \sqrt[3]{k} \, \exp\left( -i \, \pi/3\right)}.
\end{equation*}
If $k\to 0^+$, then $z_1$ and $z_2= \overline{z}_1$ converge to $z=1$ inside the unit disc.
Now, for each $\delta>0$ there exists a solution of \eqref{eq:de2} having two distinct zeros in $D(1,\delta) \cap \D$.
\end{example}

\smallskip

The following example concerns the sharpness of Theorem~\ref{thm:growth}. 


\begin{example} \label{ex:growth}
Let $0 \leq K <\infty$, and let $A$ be the analytic function 
\begin{equation*}
A(z) = - \frac{K+4\sqrt{1+K} \left( \log\frac{e}{1-z} \right)^{-1}}{4(1-z)^2}, \quad z\in\D.
\end{equation*}
Now
\begin{equation*}
\abs{A(z)} (1-\abs{z}^2)^2 
       \leq K + \frac{4\sqrt{1+K}}{\log\frac{e}{1-\abs{z}}} , \quad z \in \D.
\end{equation*}
However, the analytic function
\begin{equation*}
f(z) = \frac{1}{(1-z)^{(\sqrt{1+K}-1)/2}} \, \log\frac{e}{1-z}
\end{equation*}
is a solution of \eqref{eq:de2} such that
$f\not\in H^\infty_p$ for $p=(\sqrt{1+K}-1)/2$.
\end{example}

\smallskip

The following example shows that the number of zeros of a solution of
\eqref{eq:de2} may be larger than any pregiven number, while the
coefficient function $A$ satisfies \eqref{eq:cs} for some
sufficiently large $0<C<\infty$.  See \cite{CDO:2008} for
similar examples concerning the cases of $A\in H^\infty_0$ and $A\in H^\infty_1$.
Before the example, we recall some basic
properties of the Legendre polynomials $P_0, P_1, P_2, \dotsc$, which
 can be recovered from Bonnet's recursion formula
\begin{equation*}
n P_n (z) = (2n - 1) z P_{n-1}(z) - (n-1) P_{n-2}(z), \quad P_0(z) = 1, \quad P_1(z) = z.
\end{equation*}
For every $n\in\N$, the Legendre polynomial $P_n$ is known to have $n$ distinct zeros in the interval $(-1,1)$,
and $P_n$ is a  solution of Legendre's differential equation
\begin{equation} \label{eq:LegendreDiff}
(1-z^2) P''_n(z) -2 z P'_n(z) + n(n+1) P_n(z) =0, \quad z\in\D, \quad n\in\N \cup \{0\}.
\end{equation}


\begin{example} \label{ex:csfinite}
Let $P_n$ be the Legendre polynomial for $n\in \N\cup \{0\}$. By \eqref{eq:LegendreDiff},
\begin{equation*} 
P''_n(z) + a_1(z) P'_n(z) + a_0(z) P_n(z) =0, \quad a_1(z) = \frac{-2z}{1-z^2}, \quad a_0(z)=\frac{n(n+1)}{1-z^2}, 
\end{equation*}
for any $z\in\D$. Define $b(z) = -(1/2) \log{(1-z^2)}$ for $z\in\D$, and note that then
$b$ is a~primitive of $-a_1/2$. 
According to \cite[p.~74]{L:1993} the analytic function
\begin{equation*}
f(z) = P_n(z) e^{-b(z)} = P_n(z)\,  (1-z^2)^{1/2}, \quad z\in\D,
\end{equation*}
which is bounded and has precisely $n$ zeros in $\D$, is a solution of \eqref{eq:de2} with
\begin{equation*}
  A(z)=a_0(z) -\frac{1}{4} \, \big( a_1(z) \big)^2 - \frac{1}{2} \, a_1'(z) = \frac{1+n(n+1)(1-z^2)}{(1-z^2)^2}, 
  \quad z\in\D.
\end{equation*}

Let us consider the case $n=0$ more closely. 
It is easy to verify that a solution base $\{f_1,f_2\}$ of \eqref{eq:de2} with $A(z) = (1-z^2)^{-2}$
is given by the non-vanishing functions
\begin{equation*}
f_1(z) = (1-z^2)^{1/2}, \quad f_2(z) = (1-z^2)^{1/2} \, \log \frac{1+z}{1-z}, \quad z\in\D.
\end{equation*}
Let $0<\alpha<\infty$. Then
\begin{equation*}
f(z) = f_1(z) - \frac{1}{\alpha} \, f_2(z) = (1-z^2)^{1/2} \left( 1 - \frac{1}{\alpha} \, \log\frac{1+z}{1-z} \right)
\end{equation*}
is also a solution of \eqref{eq:de2}. Evidently, the point $z=z(\alpha)$ is a zero of $f$ 
if and only if $z(\alpha) = (e^\alpha-1)/(e^\alpha+1)$. Since
$z(\alpha)$ is a continuous function of $0<\alpha<\infty$, we conclude that there exists a
solution base $\{ f_1-f_2/\alpha_1, f_1-f_2/\alpha_2\}$ of \eqref{eq:de2}, where $0<\alpha_1<\alpha_2<\infty$,
such that the hyperbolic distance $\varrho_h(z(\alpha_1),z(\alpha_2))$ is smaller than any pregiven number.
In particular, the constant $\delta>0$ in Corollary~\ref{cor:st} depends on the choice of linearly
independent solutions,
even if $A$ satisfies \eqref{eq:cs} for $C=0$.
\end{example}

\smallskip
By considering similar examples one can investigate the sharpness of the second assertions
of Theorems~\ref{thm:separation} and \ref{th:horodiscs}. Details are left for the interested reader.

We next offer two concrete examples of equations whose solutions
admit infinite oscillation, but the coefficient satisfies
\begin{equation} \label{eq:wep}
|A(z)| (1-|z|^2)^2 \leq 1 + \varepsilon(|z|), \quad z\in\D,
\end{equation}
where $\varepsilon(|z|)$ decays to zero slower than the linear rate
as $|z|\to 1^-$. The following example is similar to \cite[Example~12]{CGHR:2012}.


\begin{example}

Let $p$ be a locally univalent analytic function in $\D$. The
functions
    $$
    f_1(z)=\big(p'(z) \big)^{-1/2} \, \sin p(z),\quad f_2(z)=\big(p'(z) \big)^{-1/2} \, \cos p(z), \quad z\in\D,
    $$
are linearly independent solutions of \eqref{eq:de2} with $A = (p')^2 + S_p/2$.
We consider the equations \eqref{eq:de2} with $A=A_1$ and $A=A_2$ induced by
    \begin{equation*}
    p_1(z)=\log\left(\log\frac{e^e}{1-z}\right),\quad
    p_2(z)=\left(\log\frac{e}{1-z}\right)^q,
    \end{equation*}
where $0<q<1$. In the first case
    \begin{equation*}
    \begin{split}
    A_1(z)(1-z)^2&=\frac14\, \frac{5+\left( \log\frac{e^e}{1-z}\right)^2}{\left( \log\frac{e^e}{1-z}\right)^2},
    \end{split}
    \end{equation*}
and it follows that \eqref{eq:wep} holds for $\e_1(r)\sim 5 (\log (e^e/(1-r)))^{-2}$ as $r\to1^-$; 
the zeros of the solution $f_1$ are $z_k = 1 - \exp( e - \exp(k\pi) )$, where $k\in\mathbb{Z}$. 
In the second case
    $$
    A_2(z)(1-z)^2=\frac14+\frac{q^2}{\left(\log\frac{e}{1-z}\right)^{2(1-q)}}
    +\frac14\, \frac{1-q^2}{\left(\log\frac{e}{1-z}\right)^2},
    $$
and so  $\e_2(r)\sim 4q^2 (\log (e/(1-r)))^{2(q-1)}$ as $r\to1^-$;
the zeros of the solution $f_1$ are $z_k=1-\exp \big( 1 -(k\pi)^{1/q} \big)$, where $k\in\mathbb{Z}$. 
\end{example}


\section{Proof of Theorem~\ref{Thm-Nehari-Local}}

The proof is grounded on an application of a suitable family of conformal maps.
The following construction, including Lemma~\ref{LemmaLocalize} below,
is borrowed from \cite[p.~576]{GGPPR2010}.
Without loss of generality, we may assume $\zeta=1$.


\begin{figure}[h!]
  \centering
  \includegraphics[width=0.4\textwidth,keepaspectratio=true]{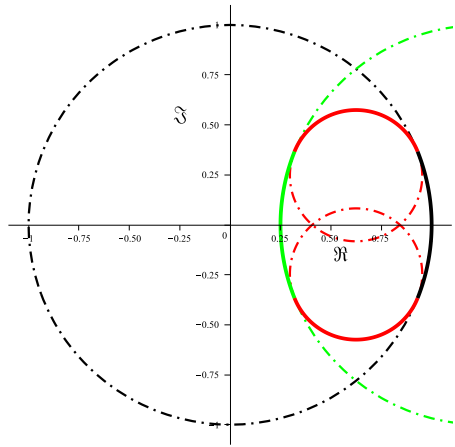}
  \caption{The boundary of $\Omega_{3/8,1/4}$ consists of the colorized bold curves.}
  \label{fig:omegadomain}
\end{figure}

Let $\tau,\r\in(0,1)$ such that $2\tau+\r<1$. Consider the circles
$\partial\D$, $\partial D(1+\r,1)$ and $\partial D(c_{\pm},r)$, where
$c_{\pm}=(1+\r)(1\pm i\tan\tau)/2$ and $r=|e^{i\tau}-c_+|$.
The discs $D(c_{\pm},r)$ are contained in both $\D$ and
$D(1+\r,1)$. Moreover, the circles $\partial D(c_{\pm},r)$
intersect $\partial\D$ on the points $e^{\pm i\tau}$, and the common
points of $\partial D(c_{\pm},r)$ and $\partial D(1+\r,1)$ are the
reflections of $e^{\pm i\tau}$ with respect to the line ${\rm Re} \, z=(1+\r)/2$. Let us call them $\gamma_\pm$ according to the
sign of their imaginary parts. Let $\Omega_{\tau,\r}$ be the
Jordan domain formed by the shortest four circular arcs connecting
$e^{\pm i\tau}$ and $\gamma_\pm$ on these four circles. See Figure~\ref{fig:omegadomain} for an illustration.
Let $\phi_{\tau,\r}$ be the conformal map of $\D$ onto
$\Omega_{\tau,\r}$. The existence of a such mapping is
ensured by the Riemann mapping theorem, which also shows
that under the additional conditions $\phi_{\tau,\r}(0)=(1+\r)/2$
and $\phi_{\tau,\r}'(0)>0$ this mapping is unique.

The following lemma produces an estimate for the growth
of the Schwarzian derivative of $\phi_{\tau,\r}$. An alternative approach
is explained in \cite[pp.~198--208]{N:1975}: the Schwarzian derivative of $\phi_{\tau,\r}$
is explicitly determined by the boundary arcs of $\Omega_{\tau,\r}$, since
the boundary $\partial\Omega_{\tau,\r}$ forms a curvilinear polygon.


\begin{letterlemma}[\protect{\cite[Lemma~8]{GGPPR2010}}]\label{LemmaLocalize}
Let $0<p<\infty$, and let $\tau,\r\in(0,1)$ such that $2\tau+\r<1$.
Then the function $\phi_{\tau,\r}$ satisfies $(\log\phi'_{\tau,\r})'\in H^p$,
$\phi''_{\tau,\r}\in H^p$ and
    \begin{equation*}
    \int_\D\bigg|\frac{\phi_{\tau,\r}''(z)}{\phi_{\tau,\r}'(z)}\bigg|^p\,dm(z)\longrightarrow 0,\quad\tau\to0^+.
    \end{equation*}
\end{letterlemma}

Before the proof of Theorem~\ref{Thm-Nehari-Local}, we make some observations about 
$\vp_{\tau,\r}$. To conclude that $\vp_{\tau,\r}$ maps the open interval $(-1,1)$ into the real axis,
we follow \cite[p.~11]{P:1975}. Evidently, $\vp_{\tau,\r}$ admits a Taylor expansion
\begin{equation*}
\vp_{\tau,\r}(z) = \frac{1+\r}{2} + a_1 z + a_2 z^2 + a_3 z^3 + \dotsb,
\end{equation*}
where $a_1>0$. Define an auxiliary function $\widetilde{\vp}(z) = \overline{\vp_{\tau,\r}(\overline{z})}$. 
Function $\widetilde{\vp}$ is analytic and univalent in $\D$, and it has a Taylor expansion
\begin{equation*}
\widetilde{\vp}(z) = \frac{1+\r}{2} + a_1 z + \overline{a}_2 z^2 +  \overline{a}_3 z^3 + \dotsb .
\end{equation*}
Note that $\widetilde{\vp}(\D) = \vp_{\tau,\r}(\D)$,
and hence $\widetilde{\vp} \equiv \vp_{\tau,\r}$ according to the uniqueness part of the Riemann mapping theorem.
Since the Taylor expansion of $\vp_{\tau,\r}$ is unique, we conclude that coefficients $a_j$ are real for all $j\in\N$.
This means that $\vp_{\tau,\r}$ maps the interval $(-1,1)$ into the real axis, and hence is typically real. 
Furthermore, since $\vp_{\tau,\r}$ is univalent, and as a real function of a real variable 
it is increasing at $z=0$ by $\vp_{\tau,\r}'(0)>0$, we have shown that $z=1$ is a fixed point of $\vp_{\tau,\r}$.
Recall that $\vp_{\tau,\r}$ has an injective and continuous extension
to the closed unit disc~$\overline{\D}$ by the famous theorem of Carath\'eodory.

Moreover, $\vp'_{\tau,\r}$ has a
continuous extension to~$\overline{\D}$, and
    \begin{equation*}
    \left|\frac{\vp_{\tau,\r}''(z)}{\vp_{\tau,\r}'(z)}\right|(1-|z|^2)^\frac1p\longrightarrow 0,\quad
    |z|\to1^-,
    \end{equation*}
since both $\vp_{\tau,\r}''$ and $(\log \vp_{\tau,\r}')'$
belong to $H^p$ for all $0<p<\infty$, see \cite[Theorems~3.11 and 5.9]{D:2000}. Standard estimates
yield
    \begin{equation}\label{Eq:SchwarzianGrowthT}
    |S_{\vp_{\tau,\r}}(z)| (1-|z|^2)^{1+\frac1p}\longrightarrow 0,\quad |z|\to1^-,
    \end{equation}
for all $1\le p<\infty$.

To prove the first assertion of Theorem~\ref{Thm-Nehari-Local},
assume that there exists $\d>0$ such that any non-trivial solution of
\eqref{eq:de2} has at most one zero in $D(1,\d)\cap\D$. Fix
$\tau,\r\in(0,1)$ such that $2\tau+\r<1$ and
$\vp_{\tau,\r}(\D)=\Omega_{\tau,\r}\subset D(1,\d) \cap \D$. Write
$T=\vp_{\tau,\r}$ for short. Then, for any given linearly
independent solutions $f_1$ and $f_2$ of \eqref{eq:de2}, the
meromorphic function $f_1/f_2\circ T$ is univalent in $\D$.
Therefore
    \begin{equation}\label{Eq:Kraus}
    \big| S_{f_1/f_2}\big(T(z)\big) \big(T'(z)\big)^2+S_T(z)\big| (1-|z|^2)^2 \leq 6,\quad z\in\D,
    \end{equation}
by Kraus theorem~\cite{Kraus1932}, see also \cite[pp.~67--68]{P:1975}
regarding the meromorphic case.

Let $\{w_n\}$ be any sequence of points in $\D$
such that $w_n\to 1$, and define $z_n$ by the equation $T(z_n)=w_n$. Then
\begin{equation*}
\lim_{n\to\infty} z_n = \lim_{n\to\infty} T^{-1}(w_n) = T^{-1}\left( \, \lim_{n\to\infty} w_n \right) = T^{-1}(1) = 1,
\end{equation*}
because $z=1$ is a fixed point of $T$. Let~$L_n$ denote
the straight line segment from $z_n\in\D$ to $z_n/\abs{z_n}\in\partial\D$. For all
$n$ sufficiently large, we have
    \begin{equation*}
    1-\abs{T(z_n)} = \left|T( z_n/\abs{z_n}) \right| -
    \abs{T(z_n)}\\
     \le \left| \int_{L_n} T'(z) \, dz \right| \leq (1-\abs{z_n}) \, \sup_{z\in L_n} \, \abs{T'(z)}.
    \end{equation*}
We point out that $T'$ is continuous in $\overline{\D}$, and $T'(1)\neq 0$
by the Julia-Carath\'eodory theorem. 
By the Schwarz-Pick lemma, we deduce
    $$
    1\leq \frac{1-|T(z_n)|^2}{|T'(z_n)|(1-|z_n|^2)}\le\frac{1+|T(z_n)|}{1+|z_n|}\, \frac{\sup_{z\in L_n} \, \abs{T'(z)}}{|T'(z_n)|}
    \longrightarrow 1,\quad
    n\to\infty.
    $$
This, together with \eqref{Eq:SchwarzianGrowthT} and
\eqref{Eq:Kraus}, yields
    \begin{equation*}
    \begin{split}
    &\limsup_{n\to\infty}\, |A(w_n)|(1-|w_n|^2)^2\\
    &\quad =\limsup_{n\to\infty}\, \frac12 \, \big| S_{f_1/f_2}(T(z_n)) \big| \big( 1-|T(z_n)|^2 \big)^2\\
    &\quad
    =\limsup_{n\to\infty}\, \frac12 \, \big|S_{f_1/f_2}(T(z_n))\big| \, |T'(z_n)|^2(1-|z_n|^2)^2\left(\frac{1-|T(z_n)|^2}{|T'(z_n)|(1-|z_n|^2)}\right)^2\le3,
    \end{split}
    \end{equation*}
which is a contradiction.


\section{Proof of Theorem~\ref{thm:growth}}

Let $f$ be a solution of \eqref{eq:de2}, and suppose that $0\leq \delta<R<1$.  We have
\begin{equation*}
  |f(z)| 
 \leq \int_\delta^{|z|} \!\!\!\int_\delta^t \, \left| f''\Big(s \, \tfrac{z}{|z|} \Big) \right| \, ds dt 
  +M(\delta,f') + M(\delta,f), \quad \delta<|z|<1,
\end{equation*}
where $M(\delta,\,\cdot\, )$ is the maximum modulus on $|z|=\delta$.
By means of \eqref{eq:de2}, we obtain
\begin{align}
& \sup_{\delta < |z| <R} \, (1-|z|^2)^p |f(z)| \notag\\
&  \qquad \leq \left( \sup_{\delta < |\zeta| <R} \, (1-|\zeta|^2)^p |f(\zeta)| \right) 
  \left( \sup_{\delta < |\zeta| <R} \, (1-|\zeta|^2)^2 |A(\zeta)| \right) \label{eq:mmest}\\
& \qquad \qquad \cdot \sup_{\delta < |z| <R} \left(  (1-|z|^2)^p \int_\delta^{|z|} \!\!\!\int_\delta^t \, \frac{ds dt}{(1-s^2)^{p+2}} \right)  + M(\delta,f') + M(\delta,f). \notag
\end{align}
If $0<p<\infty$, then
\begin{equation*}
 \lim_{|z|\to 1^-} \left(  (1-|z|^2)^p \int_\delta^{|z|} \!\!\!\int_\delta^t \, \frac{ds dt}{(1-s^2)^{p+2}} \right)
   = \frac{1}{4 p (p+1)}
\end{equation*}
by the Bernoulli-l'H\^{o}pital theorem.
The estimate \eqref{eq:mmest} implies that $\nm{f}_{H^{\infty}_p} <\infty$ provided that $K<4p(p+1)$.
We conclude $\nm{f}_{H^{\infty}_p} <\infty$ for any $p>(\sqrt{1+K}-1)/2$.


\section{Proof of Theorem~\ref{thm:separation}}

To prove the first assertion of Theorem~\ref{thm:separation},
suppose that $f_1$ is any non-trivial solution of \eqref{eq:de2} which has two distinct
zeros $z_1,z_2\in\D$ such that their hyperbolic midpoint $\xi = \xi_h(z_1,z_2)$ satisfies
$1-|\xi|<1/C$. Let $\{f_1,f_2\}$ be a solution base of \eqref{eq:de2}. Define
$h=f_1/f_2$, which implies that $S_h=2A$. 

Let $a\in\D$ such that $1-1/C<|a|<1$. If we define $r_a= 1 - C^{1/2}(1-|a|)^{1/2}$, then $0<r_a<1$.
Set $g_a(z)=(h\circ\varphi_a)(r_az)$. Then the assumption \eqref{eq:cs} yields
\begin{align*}
  |S_{g_a}(z)|(1-|z|^2)^2
  & = \big|S_h\big(\varphi_a(r_az) \big)\big| \, |\varphi'_a(r_az)|^2 \, r_a^2 \, (1-|z|^2)^2 \\
  & \le 2\Big( 1+C\big( 1-|\varphi_a(r_az)| \big) \Big) \left(\frac{1-|z|^2}{1-r_a^2|z|^2}\right)^2r_a^2 \\
  & \le 2 \Big( 1+C \big(1-\varphi_{|a|}(r_a)\big) \Big) \, r_a^2 \leq 2, \quad z\in\D,
\end{align*}
where the last inequality follows from
\begin{align*}
  1 -  \Big( 1+C \big(1-\varphi_{|a|}(r_a)\big) \Big) \, r_a^2
  & = C \, \frac{(1+r_a)(1-|a|)}{1-|a|r_a} \left( \frac{1-|a|r_a}{1-r_a} - r_a^2 \right) \geq 0.
\end{align*}
According to \cite[Theorem~1]{N:1949} the function  $g_a$ is univalent in $\D$,
and hence $h=f_1/f_2$ is univalent in the pseudo-hyperbolic disc $\Delta_p(a,r_a)$.

The  argument above shows that $h$ is univalent in $\Delta_p( \xi, r_\xi)$, and hence
\begin{equation*} 
    \varrho_h(z_1,z_2)  = \log \frac{1+\varrho_p(z_1,\xi)}{1-\varrho_p(z_1,\xi)} 
     \geq \log \frac{1+r_\xi}{1-r_\xi}
    = \log \frac{2 - C^{1/2}(1-|\xi|)^{1/2}}{C^{1/2}(1-|\xi|)^{1/2}}.
\end{equation*}

To prove the second assertion of Theorem~\ref{thm:separation},
suppose that the hyperbolic distance between any distinct zeros $z_1,z_2\in\D$ of any
non-trivial solution of \eqref{eq:de2},
for which $1-|\xi_h(z_1,z_2)|<1/C$, satisfies \eqref{eq:separation} with some $0<C<\infty$. In another words,
\begin{equation} \label{eq:iaw}
\varrho_h(z_1,z_2) \geq \log\frac{1+r_\xi}{1-r_\xi}, \quad r_\xi=1-C^{1/2} (1-|\xi|)^{1/2}, \quad \xi=\xi_h(z_1,z_2).
\end{equation}

First, we show that each non-trivial solution of \eqref{eq:de2} has at most
one zero in
\begin{equation*} 
\Delta_h\bigg( a, \,\frac{1}{2}\log\frac{1+R_a}{1-R_a}\bigg),
\quad R_a = 1 - (8C)^{1/3} (1-|a|)^{1/3}, \quad 1-|a|<(8C)^{-1}.
\end{equation*}
Assume on the contrary that there exists a non-trivial solution having two distinct
zeros $z_1,z_2 \in \Delta_p(a,R_a)$ for some $1-|a|<(8C)^{-1}$. 
By hyperbolic geometry we conclude $\xi\in\Delta_p(a,R_a)$, and hence
\begin{equation*}
1-r_\xi \leq C^{1/2}  \left( 1 - \frac{|a|-R_a}{1-|a|R_a} \right)^{1/2}
        = \frac{C^{1/2}   (1-|a|)^{1/2} \, (1+R_a)^{1/2}}{(1-|a| R_a)^{1/2}},
\end{equation*}
which implies
\begin{align*}
\frac{1+r_\xi}{1-r_\xi} \cdot \frac{1-R_a}{1+R_a}
     & \geq \frac{(1+r_\xi) (1-|a| R_a)^{1/2}}{C^{1/2}(1-|a|)^{1/2} \,(1+R_a)^{1/2}} \cdot \frac{(8C)^{1/3} (1-|a|)^{1/3}}{1+R_a}\\
     & \geq \frac{ (8C)^{1/3}}{2 \sqrt{2} \, C^{1/2}} \cdot \frac{(1- R_a)^{1/2}}{(1-|a|)^{1/6}} 
       =1 .
\end{align*}
We deduce
\begin{equation*}
\rho_h(z_1,z_2) \leq \rho_h(z_1,a) + \rho_h(a,z_2) < \log\frac{1+R_a}{1-R_a} \leq \log\frac{1+r_\xi}{1-r_\xi}, 
\end{equation*}
which contradicts \eqref{eq:iaw}.

Second, we derive the estimate \eqref{eq:hypconc}. Let $\{f_1,f_2\}$ be a solution base of \eqref{eq:de2} and set
$h=f_1/f_2$ so that $S_h=2A$. Set $g_a(z)=(h\circ\varphi_a)(R_a z)$.
Since $h$ is univalent in each pseudo-hyperbolic disc $\Delta(a,R_a)$ for $1-|a|<(8C)^{-1}$, it follows
that $g_a$ is univalent in $\D$ for those values of $a$, and hence
    \begin{equation*}
    \begin{split}
    |S_{g_a}(z)|(1-|z|^2)^2&=\big|S_h\big(\varphi_a(R_az) \big)\big| |\varphi'_a(R_a z)|^2\, R_a^2\,(1-|z|^2)^2\\
    &=2\big|A\big(\varphi_a(R_az)\big) \big| \frac{(1-\abs{a}^2)^2}{|1-\overline{a} R_a z|^4} \, R_a^2 \, (1-|z|^2)^2\le6,
    \quad z\in\D,
    \end{split}
    \end{equation*}
by Kraus theorem~\cite{Kraus1932}. By choosing $z=0$, we conclude
\begin{equation*}
    |A(a)| (1-\abs{a}^2)^2 
    \le\frac{3}{R_a^2}
    =3 \left( 1 +  (1-R_a) \, \frac{1+R_a}{R_a^2} \right),
    \quad 1-|a|<(8C)^{-1}.
\end{equation*}


\section{Proof of Theorem~\ref{th:horodiscs}}

We begin with the first assertion of Theorem~\ref{th:horodiscs}.
Suppose that the coefficient $A$ satisfies \eqref{eq:cs} for some
$0<C<\infty$, and there exists
a non-trivial solution $f_1$ of \eqref{eq:de2} having two distinct zeros $z_1,z_2 \in D(\zeta,(1+C)^{-1})$
for some $\abs{\zeta} \leq C/(1+C)$.
Let $f_2$ be a solution of \eqref{eq:de2} linearly independent to $f_1$.
By setting $h=f_1/f_2$, we deduce $S_h=2A$.
The M\"obius transformation $T(z) = \zeta + (1+C)^{-1} z$ is a conformal map
from $\D$ onto $D(\zeta,(1+C)^{-1})$.
We proceed to prove that $g = h\circ T$ is univalent in $\D$. Now
\begin{align*}
    |S_g(z)|(1-|z|^2)^2 & = \big|S_h\big(T(z) \big)\big|\, |T'(z)|^2 (1-|z|^2)^2
                           \notag\\
    & \leq 2 \, \frac{1 + C \big( 1-\abs{\zeta + (1+C)^{-1}z} \big)}{\big( 1- \abs{\zeta + (1+C)^{-1}z}^2 \big)^2}
    \,(1+C)^{-2}  (1-|z|^2)^2, \quad z\in\D.  
\end{align*}

The proof of the first assertion is divided into two separate cases.
By differentiation, there exists $0<t_C<1/3$ such that
the auxiliary function $$\mu(t) = \big( 1+C(1-t) \big) (1-t^2)^{-2}, \quad 0<t<1,$$
is decreasing for $0<t<t_C$, and increasing for $t_C<t<1$.

(i) Suppose that $z\in\D$, and $\abs{\zeta+(1+C)^{-1} z}>t_C$. By the triangle inequality
$\abs{\zeta+(1+C)^{-1}z} \leq (1+C)^{-1} (C + \abs{z})$, and hence
\begin{align}
  |S_g(z)|(1-|z|^2)^2
  & \leq 2 \, \frac{1 + C \big( 1-  (1+C)^{-1} (C + \abs{z}) \big)}{\big( 1- (1+C)^{-2} (C + \abs{z})^2 \big)^2}
  \, (1+C)^{-2}  (1-|z|^2)^2 \notag\\
  & = 2\, \frac{(1+2C-C\abs{z}) (1+C) (1+\abs{z})^2}{(\abs{z}+2C+1)^2}
  \leq  2. \label{eq:incr}
\end{align}
The inequality in \eqref{eq:incr} follows by differentiation, since the quotient 
is an increasing function of $\abs{z}$ for $0<\abs{z} < 1$.

(ii) Suppose that $z\in\D$, and $\abs{\zeta+(1-a)z}\leq t_C$. Since $\mu(0)\geq \mu(t)$
for all $0<t\leq t_C$, we deduce
\begin{equation*}
  |S_g(z)|(1-|z|^2)^2
  \leq \frac{2}{1+C} < 2.
\end{equation*}
By means of (i), (ii) and \cite[Theorem~1]{N:1949}, we conclude that $g$ is univalent in $\D$.
This is a contradiction, since the preimages $T^{-1}(z_1) \in\D$ and $T^{-1}(z_2)\in\D$ are distinct zeros of $g$.
The first assertion of Theorem~\ref{th:horodiscs} follows.

We turn to consider the second assertion of Theorem~\ref{th:horodiscs}.
Suppose that there exists $0<C<\infty$ such that any non-trivial solution of \eqref{eq:de2}
has at most one zero in any Euclidean disc $D(\zeta,(1+C)^{-1})$ for  $|\zeta| \leq C/(1+C)$.
Let $\{f_1,f_2\}$ be a solution base of \eqref{eq:de2}, and set
$h=f_1/f_2$ which implies $S_h=2A$. 
Suppose that $a\in\D$ and $|a|>C/(1+C)$.
Set $g_a(z)=(h\circ\varphi_a)(r_a z)$, where 
    \begin{equation*}
    r^2_a= \frac{\abs{a}-\frac{C}{1+C}}{\abs{a} \left( 1-\abs{a} \, \frac{C}{1+C} \right)},
    \quad 0<r_a<1.
    \end{equation*}
Then
\begin{equation*}
\Delta_p(a,r_a)
=D \bigg( \frac{a}{\abs{a}}\cdot \frac{C}{1+C}, \, \frac{r_a(1-\abs{a}^2)}{1-r_a^2 \abs{a}^2} \bigg)
\subset D\left( \frac{a}{\abs{a}}\cdot \frac{C}{1+C}, \, \frac{1}{1+C} \right).
\end{equation*}
Since $z\mapsto \varphi_a(r_a z)$ maps $\D$ onto $\Delta_p(a,r_a)$, it follows
that $g_a$ is univalent in $\D$ by the assumption. Hence
    \begin{equation*}
    \begin{split}
    |S_{g_a}(z)|(1-|z|^2)^2&=\big| S_h\big(\varphi_a(r_a z)\big) \big| \, |\varphi'_a (r_a z)|^2 \, r_a^2 \, (1-|z|^2)^2\\
    &=2 \, \big| A\big(\varphi_a(r_a z)\big) \big|  \left(\frac{1-\abs{a}^2}
      {|1-\overline{a} r_a z|^2}\right)^2 r_a^2 \, (1-|z|^2)\le6, \quad z\in\D,
    \end{split}
    \end{equation*}
by Kraus theorem~\cite{Kraus1932}. By choosing $z=0$, we obtain
    $$
    |A(a)| (1-\abs{a}^2)^2 \le\frac{3}{r_a^2}
    =3 \left( 1 +\frac{\frac{C}{1+C} \, (1+\abs{a})}{\abs{a}-\frac{C}{1+C}}  (1-\abs{a})  \right),
    \quad \frac{C}{1+C}< \abs{a} < 1.
    $$


\section{Proof of Theorem~\ref{thm:carleson}}

If $A$ is analytic in $\D$, and satisfies \eqref{eq:cs} for some $0<C<\infty$, then
any non-trivial solution $f$ of \eqref{eq:de2} has at most one zero $z_n$ in any horodisc
\begin{equation*}
\mathcal{D}_\theta = D\left( e^{i\theta} \, \frac{C}{1+C}, \, \frac{1}{1+C} \right), \quad e^{i\theta} \in \partial\D,
\end{equation*}
by Theorem~\ref{th:horodiscs}. Suppose that $0<1-r<2/(1+C)$. Now $r\in \partial\mathcal{D}_\theta$
if and only if
\begin{equation*}
\left| e^{i\theta} \, \frac{C}{1+C} - r \right|^2 = \frac{1}{(1+C)^2}.
\end{equation*}
The positive solution $\theta=\theta(r)$ of this equation satisfies
\begin{equation*}
\theta(r) = \arccos \frac{C-1 + r^2(1+C)}{2rC} \sim \sqrt{2/C} \, (1-r)^{1/2}, \quad r\to 1^{-}.
\end{equation*}
This implies that the zeros $z_n$ of $f$, for which $0<1-|z_n|<2/(1+C)$, induce
pairwise disjoint zero-free tent-like domains $\Omega_n\subset\D$, which intersect $\partial\D$ on arcs 
$I_n = \overline{\Omega}_n \cap \partial\D$ of normalized length 
\begin{equation*}
\ell( I_n) \sim \frac{\sqrt{2}}{\pi \sqrt{C}} \, (1-|z_n|)^{1/2}, \quad n\to\infty.
\end{equation*}
Consequently, if $Q$ is any Carleson square for which $\ell(Q)<2/(1+C)$, then
\begin{equation*}
\sum_{z_n\in Q} (1-|z_n|)^{1/2} \lesssim \sum_{z_n\in Q} \ell(I_n) \lesssim \ell(Q)^{1/2},
\end{equation*}
where the comparison constants depend on $0<C<\infty$. 
By a standard argument, 
this implies \eqref{eq:carleson} for any Carleson square $Q$.


\end{document}